\documentclass[12pt]{article}
\usepackage[top=1in,bottom=1in,left=1in,right=1in]{geometry}
\usepackage{indentfirst}
\usepackage{amsfonts,amsmath,amsthm,amssymb,hyperref,enumerate,bm}
\usepackage{longtable}
\usepackage[caption=false]{subfig}
\usepackage{leftidx}
\usepackage{lineno}
\usepackage{color,xcolor}
\usepackage{array}
\usepackage{amsmath}
\usepackage{latexsym}
\usepackage{float}
\usepackage{color}
\usepackage{longtable,supertabular,booktabs}
\usepackage{rotating}
\usepackage{cases}
\usepackage{multicol,multirow}
\usepackage{epsfig}
\usepackage{fancyhdr}       
\usepackage{dsfont}
\usepackage{hyperref}
\hypersetup{
	colorlinks=true,
	linkcolor=red,
	filecolor=blue,
	urlcolor=red,
	citecolor=blue,
}

\newtheorem{theorem}{Theorem}[section]

\newtheorem{lemma}{Lemma}[section]

\newtheorem{corollary}{Corollary}[section]

\newtheorem{definition}{Definition}[section]

\newtheorem{example}{Example}[section]

\def\G1{G^\mathcal{C}}

 \newenvironment{prof}{\trivlist
      \item[\hskip\labelsep
      {\itshape Proof.}]\normalfont}
      {\hspace*{\fill}$\Box$\endtrivlist}
\begin{document}
\title{Extension of the Fundamental Theorem of Algebra to Polynomial Matrix Equations over $Q$-Circulant Matrices}
\author{
Hongjian Li\footnote{E-mail\,$:$ lhj@gdufs.edu.cn. Supported by the Project of Guangdong University of Foreign Studies (Grant No. 2024RC063).}\\
{\small\it  School of Mathematics and Statistics, Guangdong University of Foreign Studies,}\\
{\small\it Guangzhou 510006, Guangdong, P. R. China} \\
}
\date{}
 \maketitle
\date{}

\noindent{\bf Abstract}\quad In this paper, we establish an analogue of the Fundamental Theorem of Algebra for polynomial matrix equations, where both the coefficient matrices and the unknown matrix are $Q$-circulant matrices. This result generalizes Abramov's result \cite{Abramov} for circulant matrices.


\medskip
\noindent{\bf MR(2020) Subject Classification} 15A24, 15B05

\section{Introduction}
The Fundamental Theorem of Algebra (or FTA) states that every nonconstant complex polynomial has at least one complex root. Some versions of the statement of the FTA first appeared early in the 17th century in the writings of several mathematicians, such as Peter Roth, Albert Girard, and Ren\'{e} Descartes. The first attempt to prove this theorem was made in 1746 by Jean Le Rond d'Alembert \cite{d'Alembert}, but his proof was not very rigorous. Attempts were also made by Leonhard Euler, Fran\c{c}ois Daviet de Foncenex, Carl Friedrich Gauss, Joseph-Louis Lagrange, Pierre-Simon de Laplace, and James Wood. However, Carl Friedrich Gauss is often credited with producing the first correct proof in his 1799 doctoral dissertation, although his proof also had some small gaps. The first textbook containing a full proof of the FTA is \textit{Cours d'analyse de l'\'{E}cole Royale Polytechnique} by  Augustin-Louis Cauchy \cite{Cauchy}. Nowadays there are many proofs of the FTA, mainly using algebra, analysis and topology (see \cite{Fine}).

Let us recall the definition of circulant matrices. A matrix $A$ is called a circulant matrix if it has the following form:
\[
A=
\begin{pmatrix}
a_0 & a_1 & a_2 &\cdots & a_{d-2} & a_{d-1} \\
a_{d-1} & a_0 & a_1 &\cdots & a_{d-3} & a_{d-2} \\
a_{d-2} & a_{d-1} & a_0 &\cdots & a_{d-4} & a_{d-3} \\
\vdots & \vdots & \vdots &\ddots & \vdots & \vdots \\
a_{2} & a_{3} & a_{4} &\cdots & a_0 & a_1 \\
a_1 & a_2 & a_3 &\cdots & a_{d-1} & a_0
\end{pmatrix},
\]
where $a_0,\,a_1,\,\dots,\,a_{d-1}\in\mathbb{C}$. It is well-known that $A$ is a circulant matrix if and only if $AC_d=C_dA$, where
\begin{equation}\label{e3}
C_d=
\begin{pmatrix}
0 & 1 & 0 &\cdots & 0 & 0\\
0 & 0 & 1 &\cdots & 0 & 0\\
\cdots& \cdots & \cdots &\cdots & \cdots & \cdots \\
0 & 0 & 0 &\cdots & 1 & 0\\
0 & 0 & 0 & \cdots & 0 & 1 \\
1 & 0 & 0 & \cdots & 0 & 0
\end{pmatrix}_{d\times d}.
\end{equation}
There have been a number of papers on circulant matrices. We refer the reader to \cite{Davis, Kalman, Kra}. Recently, the authors \cite{Li} generalized the concept of circulant matrices and defined a new class of matrices called $Q$-circulant matrices. Let $Q\in M_d(\mathbb{C})$ be a non-derogatory matrix, i.e., a matrix for which each of its eigenvalues has geometric multiplicity $1$ (\cite[Definition 3.2.4.1]{Horn}). Then a $d\times d$ complex matrix $A$ is called a $Q$-circulant matrix if $AQ=QA$ (\cite[Definition 1.1]{Li}). From a direct computation, we have that the characteristic polynomial of $C_d$ is $f(x)=x^d-1$, which has $d$ distinct roots. This means that $C_d$ is non-derogatory. So it follows from the preceding definition that  $C_d$-circulant matrices are circulant matrices. For convenience, we denote the set
\[
\{X\in M_d(\mathbb{C}): AX=XA\}
\]
by $C(A)$ throughout this paper, where $A\in M_d(\mathbb{C})$ is a given matrix. By this notation, $C(C_d)$ denotes the set of all $d\times d$ complex circulant matrices. In \cite{Li}, the authors gave the algebraic structure of $Q$-circulant matrices as follows:

\begin{lemma}{\rm (\cite[Theorem 2.1]{Li})}\label{le1}
Let $Q\in M_d(\mathbb{C})$ be a non-derogatory matrix. Then $C(Q)$ is a commutative $\mathbb{C}$-algebra, and $\left\{I,\,Q,\,Q^2,\,\dots,\,Q^{d-1}\right\}$ is a $\mathbb{C}$-basis of $C(Q)$.
\end{lemma}

By Lemma \ref{le1}, we have that $C(Q)$ is a $d$-dimensional vector space over the complex field. Moreover, $\left\{I,\,Q,\,Q^2,\,\dots,\,Q^{d-1}\right\}$ is a $\mathbb{C}$-basis of $C(Q)$. Let $A$ be a $Q$-circulant matrix. Then the representation of $A$ in the basis $\left\{I,\,Q,\,Q^2,\,\dots,\,Q^{d-1}\right\}$ is unique. It follows that every circulant  matrix  is
a linear combination of the powers of the matrix $C_d$. In \cite{Li}, the authors gave the definition of the representation polynomial as follows:
\begin{definition}{\rm (\cite[Definition 2.1]{Li})}
Let $Q\in M_d(\mathbb{C})$ be a non-derogatory matrix. Let $A$ be a $Q$-circulant matrix, and suppose that $A=\sum_{i=0}^{d-1}a_iQ^i$ for some $a_i\in\mathbb{C},\,i=0,\,1,\,\dots,\,d-1$. Then $P_A(x)=\sum_{i=0}^{d-1}a_ix^i$ is called the representation polynomial of $A$.
\end{definition}

In \cite{Abramov}, V. M. Abramov established an analogue of the FTA for polynomial matrix equations, where both the coefficient matrices and the unknown matrix are circulant matrices. In \cite{Abramov2}, he further extended certain characterization theorems from the theory of functions of a complex variable to an analogous theory of functions over the algebraic structure of circulant matrices, which generalized \cite[Theorem 2]{Abramov}. In this paper, we focus on extending the FTA to polynomial matrix equations over $Q$-circulant matrices. Specifically, we present the following theorem.

%

\begin{theorem}\label{th1}
Let $Q$ be a $d\times d$ complex matrix which has $d$ distinct eigenvalues $\lambda_1,\lambda_2,\dots,\lambda_d$. Suppose that $T^{-1}QT={\rm diag}\left(\lambda_1,\,\lambda_2,\,\dots,\,\lambda_d\right)$. Then all solutions of the equation
\begin{equation}\label{e1}
X^n+A_1X^{n-1}+A_2X^{n-2}+\dots+A_{n-1}X+A_{n}=O
\end{equation}
in $C(Q)$ are given by $\left\{T{\rm diag}\left(u_1,\,u_2,\,\dots,\,u_d\right)T^{-1}: g_i(u_i)=0,i=1,2,\dots,d\right\}$, where
\begin{equation}\label{e2}
g_i(x)=x^n+f_1(\lambda_i)x^{n-1}+f_2(\lambda_i)x^{n-2}+\dots+f_{n-1}(\lambda_i)x+f_{n}(\lambda_i),
\end{equation}
and $f_1(x),f_2(x),\dots,f_n(x)$ are the representation polynomials of $A_1,A_2,\dots,A_n$, respectively.
\end{theorem}

By Theorem \ref{th1}, we obtain the following corollary.

\begin{corollary}
Let $n_i$ denote the number of distinct roots of equation \eqref{e2}. Then the total number of solutions of equation \eqref{e1} in $C(Q)$ is $\prod_{i=1}^{d} n_i$.
\end{corollary}

\section{Proof of Theorem \ref{th1}}
In this section, we give the proof of Theorem \ref{th1}. Before proving Theorem \ref{th1}, we present the following lemma.
\begin{lemma}\label{le2}
Let $Q$ be a $d\times d$ complex matrix which has $d$ distinct eigenvalues $\lambda_1,\lambda_2,\dots,\lambda_d$. Suppose that $T^{-1}QT={\rm diag}\left(\lambda_1,\,\lambda_2,\,\dots,\,\lambda_d\right)$. Then the following statements hold.
\begin{enumerate}
\item[\rm(i)]\, For any $A\in C(Q)$, we have $T^{-1}AT={\rm diag}\left(f_A(\lambda_1),\,f_A(\lambda_2),\,\dots,\,f_A(\lambda_d)\right)$, where $f_A(x)$ is the representation polynomial of $A$.
\item[\rm(ii)]\, $C(Q)=\left\{T{\rm diag}\left(u_1,\,u_2,\,\dots,\,u_d\right)T^{-1}: u_i\in\mathbb{C}, i=1,2,\dots,d\right\}$.
\end{enumerate}
\end{lemma}

\begin{prof}
We now prove the assertion (i). Since $A\in C(Q)$, it follows from Lemma \ref{le1} that $A=\sum_{i=0}^{d-1}a_iQ^i$ for some $a_i\in\mathbb{C},\,i=0,\,1,\,\dots,\,d-1$. It follows that the representation polynomial of $A$ is $f_A(x)=\sum_{i=0}^{d-1}a_ix^i$. So we have
\[
\begin{aligned}
T^{-1}AT
&=T^{-1}\left(\sum_{i=0}^{d-1}a_iQ^i\right)T
=\sum_{i=0}^{d-1}a_iT^{-1}Q^iT\\
&=\sum_{i=0}^{d-1}a_i\left(T^{-1}QT\right)^i=\sum_{i=0}^{d-1}a_i{\rm diag}\left(\lambda_1^i,\,\lambda_2^i,\,\dots,\,\lambda_d^i\right)\\
&={\rm diag}\left(\sum_{i=0}^{d-1}a_i\lambda_1^i,\,\sum_{i=0}^{d-1}a_i\lambda_2^i,\,\dots,\,
    \sum_{i=0}^{d-1}a_i\lambda_d^i\right)\\
&={\rm diag}\left(f_A\left(\lambda_1\right),f_A\left(\lambda_2\right),\dots,\,f_A\left(\lambda_d\right)\right).
\end{aligned}
\]

We next prove the assertion (ii). Let $X\in C(Q)$. Then it follows from (i) that
\[
T^{-1}XT={\rm diag}\left(f_X(\lambda_1),\,f_X(\lambda_2),\,\dots,\,f_X(\lambda_d)\right),
\]
where $f_X(x)$ is the representation polynomial of $X$. So we have
\[
X=T{\rm diag}\left(f_X(\lambda_1),\,f_X(\lambda_2),\,\dots,\,f_X(\lambda_d)\right)T^{-1},
\]
which implies that $X\in\left\{T{\rm diag}\left(u_1,\,u_2,\,\dots,\,u_d\right)T^{-1}: u_i\in\mathbb{C}, i=1,2,\dots,d\right\}$. Hence, we have
\begin{equation}\label{e10}
C(Q)\subseteq\left\{T{\rm diag}\left(u_1,\,u_2,\,\dots,\,u_d\right)T^{-1}: u_i\in\mathbb{C}, i=1,2,\dots,d\right\}.
\end{equation}
On the other hand, if $Y\in\left\{T{\rm diag}\left(u_1,\,u_2,\,\dots,\,u_d\right)T^{-1}: u_i\in\mathbb{C}, i=1,2,\dots,d\right\}$, then
\begin{equation}\label{e7}
Y=T{\rm diag}\left(u_1,\,u_2,\,\dots,\,u_d\right)T^{-1}
\end{equation}
for some $u_i\in\mathbb{C}, i=1,2,\dots,d$. Since $T^{-1}QT={\rm diag}\left(\lambda_1,\,\lambda_2,\,\dots,\,\lambda_d\right)$, we have
\begin{equation}\label{e8}
Q=T{\rm diag}\left(\lambda_1,\,\lambda_2,\,\dots,\,\lambda_d\right)T^{-1}.
\end{equation}
By \eqref{e7} and \eqref{e8}, we obtain
\[
\begin{aligned}
YQ
&=\left(T{\rm diag}\left(u_1,\,u_2,\,\dots,\,u_d\right)T^{-1}\right)\cdot\left(T{\rm diag}\left(\lambda_1,\,\lambda_2,\,\dots,\,\lambda_d\right)T^{-1}\right)\\
&=T{\rm diag}\left(u_1,\,u_2,\,\dots,\,u_d\right)\cdot{\rm diag}\left(\lambda_1,\,\lambda_2,\,\dots,\,\lambda_d\right)T^{-1}\\
&=T{\rm diag}\left(u_1\lambda_1,\,u_2\lambda_2,\,\dots,\,u_d\lambda_d\right)T^{-1}\\
&=T{\rm diag}\left(\lambda_1,\,\lambda_2,\,\dots,\,\lambda_d\right)\cdot{\rm diag}\left(u_1,\,u_2,\,\dots,\,u_d\right)T^{-1}\\
&=\left(T{\rm diag}\left(\lambda_1,\,\lambda_2,\,\dots,\,\lambda_d\right)T^{-1}\right)\cdot\left(T{\rm diag}\left(u_1,\,u_2,\,\dots,\,u_d\right)T^{-1}\right)\\
&=QY.
\end{aligned}
\]
This  implies that $Y\in C(Q)$. It follows that
\begin{equation}\label{e9}
\left\{T{\rm diag}\left(u_1,\,u_2,\,\dots,\,u_d\right)T^{-1}: u_i\in\mathbb{C}, i=1,2,\dots,d\right\}\subseteq C(Q).
\end{equation}
By \eqref{e10} and \eqref{e9}, we have $C(Q)=\left\{T{\rm diag}\left(u_1,\,u_2,\,\dots,\,u_d\right)T^{-1}: u_i\in\mathbb{C}, i=1,2,\dots,d\right\}$.

\end{prof}

We now present the proof of Theorem \ref{th1}.
\begin{prof}
Let $X\in C(Q)$. Then it follows from Lemma \ref{le2} (ii) that $X=T{\rm diag}\left(u_1,\,u_2,\,\dots,\,u_d\right)T^{-1}$ for some $u_i\in\mathbb{C},i=1,2,\dots,d$. Then we have
\begin{equation}\label{e4}
T^{-1}XT={\rm diag}\left(u_1,\,u_2,\,\dots,\,u_d\right).
\end{equation}
Let $f_1(x),f_2(x),\dots,f_n(x)$ be the representation polynomials of $A_1,A_2,\dots,A_n$, respectively. Then it follows from Lemma \ref{le2} (i) that
\begin{equation}\label{e5}
T^{-1}A_kT={\rm diag}\left(f_k(\lambda_1),\,f_k(\lambda_2),\,\dots,\,f_k(\lambda_d)\right),\quad k=1,2,\dots,n.
\end{equation}
By \eqref{e4} and \eqref{e5}, we have
\begin{equation}\label{e6}
\begin{aligned}
T^{-1}&\left(X^n+A_1X^{n-1}+A_2X^{n-2}+\dots+A_{n-1}X+A_{n}\right)T\\
&=T^{-1}X^nT+T^{-1}A_1X^{n-1}T+T^{-1}A_2X^{n-2}T+\dots+T^{-1}A_{n-1}XT+T^{-1}A_{n}T\\
&=T^{-1}X^nT+\left(T^{-1}A_1T\right)\left(T^{-1}X^{n-1}T\right)+\left(T^{-1}A_2T\right)\left(T^{-1}X^{n-2}T\right)\\
&\quad +\dots+\left(T^{-1}A_{n-1}T\right)\left(T^{-1}XT\right)+T^{-1}A_{n}T\\
&=\left(T^{-1}XT\right)^n+\left(T^{-1}A_1T\right)\left(T^{-1}XT\right)^{n-1}+\left(T^{-1}A_2T\right)\left(T^{-1}XT\right)^{n-2}\\
&\quad +\dots+\left(T^{-1}A_{n-1}T\right)\left(T^{-1}XT\right)+T^{-1}A_{n}T\\
&={\rm diag}\left(g_1(u_1),\,g_2(u_2),\,\dots,\,g_d(u_d)\right),
\end{aligned}
\end{equation}
where $g_i(x)=x^n+f_1(\lambda_i)x^{n-1}+f_2(\lambda_i)x^{n-2}+\dots+f_{n-1}(\lambda_i)x+f_{n}(\lambda_i),i=1,2,\dots,d$. It follows from \eqref{e6} that $X$ is a solution of equation \eqref{e1} if and only if $g_i(u_i)=0$ for $i=1,2,\dots,d$. Hence, all solutions of equation \eqref{e1} in $C(Q)$ are given by
\[
\left\{T{\rm diag}\left(u_1,\,u_2,\,\dots,\,u_d\right)T^{-1}: g_i(u_i)=0,i=1,2,\dots,d\right\}.
\]
\end{prof}

\section{Some Corollaries and Examples}

In this section, we give some corollaries and examples.

\subsection{All solutions of equation \eqref{e1} over $C(Q(k_1,k_2,\dots,k_d))$}

Let
\begin{equation}\label{eq19}
Q=Q(k_1,k_2,\dots,k_d)=
\begin{pmatrix}
0 & k_1 & 0 &\cdots & 0 \\
0 & 0 & k_2 &\ddots & \vdots \\
0 & 0 & 0 &\ddots & 0 \\
\vdots & \ddots & \ddots &\ddots & k_{d-1} \\
k_d & \cdots & 0 & 0 & 0
\end{pmatrix},
\end{equation}
where $k_1,\,k_2,\,\dots,\,k_d\in\mathbb{C}\backslash\{0\}$. By a direct computation, we have that the characteristic polynomial of $Q$ is
\[
f(x)=x^d-\prod_{i=1}^dk_i,
\]
which has $d$ distinct roots. Throughout this section, let $\omega=e^{i\frac{2\pi}{d}}$ denote the primitive $d$-th root of unity. Let $k=\prod_{i=1}^dk_i$ and let $\lambda$ be any $d$-th root of $k$. Then $Q$ has $d$ distinct eigenvalues $\lambda,\lambda\omega,\lambda\omega^2,\dots,\lambda\omega^{d-1}$. We define the diagonal matrix
\begin{equation}\label{eq8}
\Lambda=\frac{k}{k_d}{\rm diag}\left(1,\,\frac{\lambda}{k_1},\,\frac{\lambda^2}{k_1k_2},\,\dots,\,\frac{\lambda^{d-2}}{k_1\cdots k_{d-2}},\,\frac{\lambda^{d-1}}{k_1\cdots k_{d-1}}\right).
\end{equation}
Let $F$ be the discrete Fourier matrix, that is,
\begin{equation}\label{eq9}
F=
\frac{1}{\sqrt{d}}\begin{pmatrix}
1 & 1 & 1 &\cdots & 1 & 1\\
1 & \omega & \omega^2 &\cdots & \omega^{d-2} & \omega^{d-1}\\
1 & \omega^2 & \omega^4 &\cdots & \omega^{2\left(d-2\right)} & \omega^{2\left(d-1\right)}\\
\cdots& \cdots & \cdots &\cdots & \cdots & \cdots \\
1 & \omega^{d-2} & \omega^{2\left(d-2\right)} &\cdots & \omega^{\left(d-2\right)^2} & \omega^{\left(d-2\right)\left(d-1\right)}\\
1 & \omega^{d-1} & \omega^{2\left(d-1\right)} &\cdots & \omega^{\left(d-2\right)\left(d-1\right)} & \omega^{\left(d-1\right)^2}
\end{pmatrix}.
\end{equation}
Let $T=\Lambda F^{-1}$. Then we have
\begin{equation}\label{e11}
\begin{aligned}
T^{-1}QT
&=F\Lambda^{-1}Q\Lambda F^{-1}=F\left(\Lambda^{-1}Q\Lambda\right)F^{-1}=F\left(\lambda C_d\right)F^{-1}=\lambda FC_dF^{-1}\\
&=\lambda{\rm diag}\left(\omega^{d},\,\omega^{d-1},\,\dots,\,\omega^{2},\,\omega^{1}\right)={\rm diag}\left(\lambda\omega^{d},\,\lambda\omega^{d-1},\,\dots,\,\lambda\omega^{2},\,\lambda\omega^{1}\right),
\end{aligned}
\end{equation}
where $C_d$ is the matrix as in \eqref{e3}. By Theorem \ref{th1} and \eqref{e11}, we have the following corollary.
\begin{corollary}\label{co1}
Let $Q$ be the matrix as in \eqref{eq19}. Let $\Lambda$ and $F$ be the matrices as in \eqref{eq8} and \eqref{eq9}, respectively, and let $T=\Lambda F^{-1}$. Then all solutions of equation \eqref{e1} in $C(Q)$ are given by $\left\{T{\rm diag}\left(u_1,\,u_2,\,\dots,\,u_d\right)T^{-1}: g_i(u_i)=0,i=1,2,\dots,d\right\}$, where
\begin{equation}\label{eq1}
g_i(x)=x^n+f_1\left(\lambda\omega^{d-i+1}\right)x^{n-1}+f_2\left(\lambda\omega^{d-i+1}\right)x^{n-2}+\dots+f_{n-1}\left(\lambda\omega^{d-i+1}\right)x+f_{n}\left(\lambda\omega^{d-i+1}\right),
\end{equation}
and $f_1(x),f_2(x),\dots,f_n(x)$ are the representation polynomials of $A_1,A_2,\dots,A_n$, respectively.
\end{corollary}

Let $k_1=k_2=\dots=k_d=1$ in \eqref{eq19}. Then $C(Q)$ is the set of all $d\times d$ complex circulant matrices. By Corollary \ref{co1}, we have the following corollary.
\begin{corollary}\label{co2}
Let $F$ be the matrix as in \eqref{eq9}. Let $A_1,A_2,\dots,A_n$ be circulant matrices, and suppose that the representation polynomial of $A_k$ is $f_k(x)=\sum_{j=1}^{d}a_{k,j-1}x^{j-1}$ for $k=1,2,\dots,n$. Then all solutions of equation \eqref{e1} in circulant matrices are given by
\[
\left\{F^{-1}{\rm diag}\left(u_1,\,u_2,\,\dots,\,u_d\right)F: g_i(u_i)=0,i=1,2,\dots,d\right\},
\]
where $g_i(x)=x^n+b_1^{(i)}x^{n-1}+b_2^{(i)}x^{n-2}+\dots+b_{n-1}^{(i)}x+b_{n}^{(i)}$ and $b_k^{(i)}=\sum_{j=1}^{d}a_{k,j-1}\frac{1}{\omega^{(i-1)(j-1)}}$.
\end{corollary}

\begin{prof}
Let $k_1=k_2=\dots=k_d=1$. Then $k=1$ and $\lambda=1$. So the matrix $\Lambda=I$ becomes the identity matrix. Then $T=\Lambda F^{-1}=F^{-1}$. Note that
\begin{equation}\label{eq2}
\begin{aligned}
b_k^{(i)}&=\sum_{j=1}^{d}a_{k,j-1}\frac{1}{\omega^{(i-1)(j-1)}}=\sum_{j=1}^{d}a_{k,j-1}\omega^{(-i+1)(j-1)}\\
&=\sum_{j=1}^{d}a_{k,j-1}\left(\omega^{-i+1}\right)^{j-1}=f_k(\omega^{-i+1})=f_k(\omega^{d-i+1}).
\end{aligned}
\end{equation}
Since the representation polynomial of $A_k$ is $f_k(x)=\sum_{j=1}^{d}a_{k,j-1}x^{j-1}$, it follows from \eqref{eq1} and \eqref{eq2} that
\begin{align*}
g_i(x)
&=x^n+f_1\left(\omega^{d-i+1}\right)x^{n-1}+f_2\left(\omega^{d-i+1}\right)x^{n-2}+\dots+f_{n-1}\left(\omega^{d-i+1}\right)x+f_{n}\left(\omega^{d-i+1}\right)\\
&=x^n+b_1^{(i)}x^{n-1}+b_2^{(i)}x^{n-2}+\dots+b_{n-1}^{(i)}x+b_{n}^{(i)}.
\end{align*}
By Corollary \ref{co1}, all solutions of equation \eqref{e1} in circulant matrices are given by
\[
\left\{F^{-1}{\rm diag}\left(u_1,\,u_2,\,\dots,\,u_d\right)F: g_i(u_i)=0,i=1,2,\dots,d\right\},
\]
where $g_i(x)=x^n+b_1^{(i)}x^{n-1}+b_2^{(i)}x^{n-2}+\dots+b_{n-1}^{(i)}x+b_{n}^{(i)}$ and $b_k^{(i)}=\sum_{j=1}^{d}a_{k,j-1}\frac{1}{\omega^{(i-1)(j-1)}}$.
\end{prof}

By Corollary \ref{co2}, we have the following corollaries, which are given by \cite{Abramov}.

\begin{corollary}{\rm (\cite[Theorem 2]{Abramov})}
Equation \eqref{e1} over circulant matrices satisfies the FTA, and the total number of its solutions does not exceed $n^d$.
\end{corollary}

\begin{corollary}{\rm (\cite[Theorem 3]{Abramov})}
Let $n_i$ denote the number of distinct roots of the monic polynomial equation
\[
x^n+b_1^{(i)}x^{n-1}+\dots+b_{n}^{(i)}=0,\quad i=1,2,\dots,d,
\]
where
\[
b_k^{(i)}=\sum_{j=1}^{d} a_{k,j-1} \overline{r}_{(i-1)(j-1)}.
\]
Then the total number of solutions of \eqref{e1} in circulant matrices is $\prod_{i=1}^{d} n_i$.
\end{corollary}

We now present an example to illustrate the results of this section.

\begin{example}
\rm
Let
\[
Q=\begin{pmatrix}
0 & 1 & 0  \\
0 & 0 & 1 \\
8 & 0 & 0 \\
\end{pmatrix}.
\]
Let $\omega=e^{i\frac{2\pi}{3}}=\frac{-1+\sqrt{-3}}{2}$ be the primitive $3$-th root of unity, and let
\[
D=\begin{pmatrix}
1 & 2\omega^3 & \left(2\omega^3\right)^2  \\
1 & 2\omega^2 & \left(2\omega^2\right)^2 \\
1 & 2\omega   & \left(2\omega\right)^2 \\
\end{pmatrix}.
\]
Let $A=a_0I+a_1Q+a_2Q^2$ and $B=b_0I+b_1Q+b_2Q^2$, where
\begin{equation}\label{eq3}
\left(a_0,a_1,a_2\right)^T=D^{-1}\left(-5,2,-3\right)^T\quad \text{and}\quad \left(b_0,b_1,b_2\right)^T=D^{-1}\left(4,1,2\right)^T.
\end{equation}
We next find the solutions of the equation
\begin{equation}\label{eq4}
X^2+AX+B=O
\end{equation}
in $C(Q)$. Note that the representation polynomials of $A$ and $B$ are $f_1(x)=a_0+a_1x+a_2x^2$ and $f_2(x)=b_0+b_1x+b_2x^2$, respectively. It follows from \eqref{eq3} that $f_1\left(2\omega^3\right)=-5,f_1\left(2\omega^2\right)=2,f_1\left(2\omega\right)=-3,f_2\left(2\omega^3\right)=4,f_2\left(2\omega^2\right)=1,f_2\left(2\omega\right)=2$. So we have
\[
\begin{aligned}
g_1(x)&=x^2+f_1\left(2\omega^3\right)x+f_2\left(2\omega^3\right)=x^2-5x+4=\left(x-1\right)\left(x-4\right),\\
g_2(x)&=x^2+f_1\left(2\omega^2\right)x+f_2\left(2\omega^2\right)=x^2+2x+1=\left(x+1\right)^2,\\
g_3(x)&=x^2+f_1\left(2\omega\right)x+f_2\left(2\omega\right)=x^2-3x+2=\left(x-1\right)\left(x-2\right).
\end{aligned}
\]
Let
\[
F=\frac{1}{\sqrt{3}}\begin{pmatrix}
1 & 1 & 1  \\
1 & \omega & \omega^2 \\
1 & \omega^2 & \omega^4 \\
\end{pmatrix},
\]
and let $T={\rm diag}\left(1,\,2,\,4\right)F^{-1}$. Then it follows from Corollary \ref{co1} that the four solutions of equation \eqref{eq4} in $C(Q)$ are
\[
\begin{aligned}
T{\rm diag}\left(1,\,-1,\,1\right)T^{-1}
&=\frac{1}{6}\begin{pmatrix}
2 & -2\omega & -\omega^2  \\
-8\omega^2 & 2 & -2\omega \\
-16\omega & -8\omega^2 & 2 \\
\end{pmatrix},\\
T{\rm diag}\left(1,\,-1,\,2\right)T^{-1}
&=\frac{1}{12}\begin{pmatrix}
8 & 4+6\omega^2 & -1-3\omega^2  \\
-8-24\omega^2 & 8 & 4+6\omega^2 \\
32+48\omega^2 & -8-24\omega^2 & 8 \\
\end{pmatrix},\\
T{\rm diag}\left(4,\,-1,\,1\right)T^{-1}
&=\frac{1}{12}\begin{pmatrix}
16 & 10+4\omega^2 & 3-2\omega^2  \\
24-16\omega^2 & 16 & 10+4\omega^2 \\
80+32\omega^2 & 24-16\omega^2 & 16 \\
\end{pmatrix},\\
T{\rm diag}\left(4,\,-1,\,2\right)T^{-1}
&=\frac{1}{12}\begin{pmatrix}
20 & 10+6\omega^2 & 2-3\omega^2  \\
16-24\omega^2 & 20 & 10+6\omega^2 \\
80+48\omega^2 & 16-24\omega^2 & 20 \\
\end{pmatrix}.
\end{aligned}
\]

\end{example}

\subsection{All solutions of equation \eqref{e1} over $C(\pi)$}

Let $\lambda_1,\lambda_2,\dots,\lambda_d$ be $d$ distinct complex numbers, and let
\[
f(x)=\left(x-\lambda_1\right)\left(x-\lambda_2\right)\cdots\left(x-\lambda_d\right)=x^d+a_{d-1}x^{d-1}+\dots+a_1x+a_0.
\]
Let
\begin{equation}\label{eq5}
\pi=
\begin{pmatrix}
0 & 1 & 0 &\cdots & 0 & 0\\
0 & 0 & 1 &\cdots & 0 & 0\\
\cdots& \cdots & \cdots &\cdots & \cdots & \cdots \\
0 & 0 & 0 &\cdots & 1 & 0\\
0 & 0 & 0 & \cdots & 0 & 1 \\
-a_0 & -a_1 & -a_2 & \cdots &-a_{d-2} & -a_{d-1}
\end{pmatrix}.
\end{equation}
Note that the matrix $\pi$ is the companion matrix of the polynomial $f(x)$. So the characteristic polynomial of $\pi$ is $f(x)$. It follows that $\pi$ has $d$ distinct eigenvalues $\lambda_1,\lambda_2,\dots,\lambda_d$. Let
\begin{equation}\label{eq6}
T=
\begin{pmatrix}
1 & 1 & 1 &\cdots & 1 & 1\\
\lambda_1 & \lambda_2 & \lambda_3 &\cdots & \lambda_{d-1} & \lambda_{d}\\
\lambda_1^2 & \lambda_2^2 & \lambda_3^2 &\cdots & \lambda_{d-1}^2 & \lambda_{d}^2\\
\cdots& \cdots & \cdots &\cdots & \cdots & \cdots \\
\lambda_1^{d-2} & \lambda_2^{d-2} & \lambda_3^{d-2} &\cdots & \lambda_{d-1}^{d-2} & \lambda_{d}^{d-2}\\
\lambda_1^{d-1} & \lambda_2^{d-1} & \lambda_3^{d-1} &\cdots & \lambda_{d-1}^{d-1} & \lambda_{d}^{d-1}
\end{pmatrix}.
\end{equation}
Then we have
\begin{equation}\label{eq0}
T^{-1}\pi T={\rm diag}\left(\lambda_1,\,\lambda_2,\,\dots,\,\lambda_{d}\right).
\end{equation}
By Theorem \ref{th1} and \eqref{eq0}, we have the following corollary.

\begin{corollary}\label{co3}
Let $\pi$ and $T$ be the matrices as in \eqref{eq5} and \eqref{eq6}, respectively. Then all solutions of equation \eqref{e1} in $C(\pi)$ are given by
\[
\left\{T{\rm diag}\left(u_1,\,u_2,\,\dots,\,u_d\right)T^{-1}: g_i(u_i)=0,i=1,2,\dots,d\right\},
\]
where $g_i(x)=x^n+f_1(\lambda_i)x^{n-1}+f_2(\lambda_i)x^{n-2}+\dots+f_{n-1}(\lambda_i)x+f_{n}(\lambda_i)$, and $f_1(x),f_2(x),\dots,f_n(x)$ are the representation polynomials of $A_1,A_2,\dots,A_n$, respectively.
\end{corollary}

We now present an example to illustrate the results of this section.

\begin{example}
\rm
Let $\lambda_1=1,\lambda_2=2,\lambda_3=3$, and let
\[
f(x)=\left(x-1\right)\left(x-2\right)\left(x-3\right)=x^3-6x^{2}+11x-6.
\]
Let
\[
\pi=\begin{pmatrix}
0 & 1 & 0  \\
0 & 0 & 1 \\
6 & -11 & 6 \\
\end{pmatrix}
\quad \text{and} \quad
D=\begin{pmatrix}
1 & 1 & 1^2  \\
1 & 2 & 2^2 \\
1 & 3   & 3^2 \\
\end{pmatrix}.
\]
Let $A=a_0I+a_1\pi+a_2\pi^2$ and $B=b_0I+b_1\pi+b_2\pi^2$, where
\begin{equation}\label{eq7}
\left(a_0,a_1,a_2\right)^T=D^{-1}\left(-5,2,-3\right)^T\quad \text{and}\quad \left(b_0,b_1,b_2\right)^T=D^{-1}\left(4,1,2\right)^T.
\end{equation}
We next find the solutions of the equation
\begin{equation}\label{eq10}
X^2+AX+B=O
\end{equation}
in $C(\pi)$. Note that the representation polynomials of $A$ and $B$ are $f_1(x)=a_0+a_1x+a_2x^2$ and $f_2(x)=b_0+b_1x+b_2x^2$, respectively. It follows from \eqref{eq7} that $f_1\left(1\right)=-5,f_1\left(2\right)=2,f_1\left(3\right)=-3,f_2\left(1\right)=4,f_2\left(2\right)=1,f_2\left(3\right)=2$. So we have
\[
\begin{aligned}
g_1(x)&=x^2+f_1\left(1\right)x+f_2\left(1\right)=x^2-5x+4=\left(x-1\right)\left(x-4\right),\\
g_2(x)&=x^2+f_1\left(2\right)x+f_2\left(2\right)=x^2+2x+1=\left(x+1\right)^2,\\
g_3(x)&=x^2+f_1\left(3\right)x+f_2\left(3\right)=x^2-3x+2=\left(x-1\right)\left(x-2\right).
\end{aligned}
\]
Let
\[
T=
\begin{pmatrix}
1 & 1 & 1  \\
1 & 2 & 3 \\
1 & 4 & 9 \\
\end{pmatrix}.
\]
Then it follows from Corollary \ref{co3} that the four solutions of equation \eqref{eq10} in $C(\pi)$ are
\[
\begin{aligned}
T{\rm diag}\left(1,\,-1,\,1\right)T^{-1}
&=\begin{pmatrix}
7 & -8 & 2  \\
12 & -15 & 4 \\
24 & -32 & 9 \\
\end{pmatrix},\\
T{\rm diag}\left(1,\,-1,\,2\right)T^{-1}
&=\frac{1}{2}\begin{pmatrix}
16 & -19 & 5  \\
30 & -39 & 11 \\
66 & -91 & 27 \\
\end{pmatrix},\\
T{\rm diag}\left(4,\,-1,\,1\right)T^{-1}
&=\frac{1}{2}\begin{pmatrix}
32 & -31 & 7  \\
42 & -45 & 11 \\
66 & -79 & 21 \\
\end{pmatrix},\\
T{\rm diag}\left(4,\,-1,\,2\right)T^{-1}
&=\begin{pmatrix}
17 & -17 & 4  \\
24 & -27 & 7 \\
42 & -53 & 15 \\
\end{pmatrix}.
\end{aligned}
\]

\end{example}

\section{Conclusion}

Let $Q$ be a $d\times d$ complex matrix which has $d$ distinct eigenvalues. In this paper, we establish an analogue of the FTA for polynomial matrix equations over $Q$-circulant matrices. That is, we give all solutions of the equation
\[
X^n+A_1X^{n-1}+A_2X^{n-2}+\dots+A_{n-1}X+A_{n}=O
\]
in $C(Q)$. This result generalizes Abramov's result \cite{Abramov} for circulant matrices.

\end{document}